\documentclass[11pt]{article}
\usepackage{color}

\usepackage[utf8]{inputenc}
\usepackage[T1]{fontenc}
\usepackage{amsmath,amssymb,epsfig,bbm}
\usepackage{comment}
\usepackage[colorlinks]{hyperref}
\usepackage{todonotes}

\newtheorem{defi}{Definition}
\newtheorem{prop}[defi]{Proposition}
\newtheorem{theo}[defi]{Theorem}
\newtheorem{conj}[defi]{Conjecture}
\newtheorem{lemm}[defi]{Lemma}
\newtheorem{coro}[defi]{Corollary}
\newtheorem{rema}[defi]{Remark}
\newtheorem{exem}[defi]{Example}
\newtheorem{exems}[defi]{Examples}

\newcommand{\bdefi}{\begin{defi}}
\newcommand{\edefi}{\end{defi}}
\newcommand{\bprop}{\begin{prop}}
\newcommand{\eprop}{\end{prop}}
\newcommand{\btheo}{\begin{theo}}
\newcommand{\etheo}{\end{theo}}
\newcommand{\blemm}{\begin{lemm}}
\newcommand{\brema}{\begin{rema}}
\newcommand{\erema}{\end{rema}}
\newcommand{\bexer}{\begin{exem}}
\newcommand{\eexer}{\end{exem}}
\newcommand{\bexems}{\begin{exems}}
\newcommand{\eexems}{\end{exems}}
\newcommand{\bconj}{\begin{conj}}
\newcommand{\econj}{\end{conj}}
\newcommand{\elemm}{\end{lemm}}
\newcommand{\bcoro}{\begin{coro}}
\newcommand{\ecoro}{\end{coro}}
\newcommand{\dem}{\noindent{\bf Proof. }}

\usepackage{mathrsfs}
\renewcommand\mathcal{\mathscr}

\newcommand{\G}{{\cal G}}
\renewcommand{\H}{{\cal H}}

\newcommand{\M}{{\cal M}}

\newcommand{\OOOO}{{\cal O}}

\newcommand{\maths}[1]{{\mathbb #1}}  

\newcommand{\RR}{\maths{R}}
\newcommand{\NN}{\maths{N}}

\newcommand{\QQ}{\maths{Q}}

\newcommand{\FF}{\maths{F}}

\newcommand{\ZZ}{\maths{Z}}
\newcommand{\PP}{\maths{P}}

\newcommand{\TT}{\maths{T}}

\newcommand{\uA}{\underline{A}}

\newcommand{\aaa}{{\mathfrak a}}

\newcommand{\ra}{\rightarrow}
\newcommand{\bs}{\backslash}
\newcommand{\ov}[1]{\overline{#1}} 
\newcommand{\wt}[1]{{\widetilde{#1}}}
\newcommand{\wh}[1]{{\widehat{#1}}}
\newcommand{\ga}{\gamma}
\newcommand{\Ga}{\Gamma}

\newcommand{\cqfd}{\hfill$\Box$}

\newcommand{\id}{\operatorname{id}}
\newcommand{\pr}{\operatorname{pr}}

\newcommand{\GL}{\operatorname{GL}}

\newcommand{\PGL}{\operatorname{PGL}}

\newcommand{\heit}{\operatorname{ht}}

\newcommand{\tr}{\operatorname{tr}}

\newcounter{fig}

\title{On continued fraction expansions \\ 
of quadratic irrationals in positive characteristic}
\author{Frédéric Paulin and Uri Shapira}

\begin{document}
\bibliographystyle{../alphanum}
\maketitle
\begin{abstract}
  Let $R=\FF_q[Y]$ be the ring of polynomials over a finite field
  $\FF_q$, let $\wh K=\FF_q((Y^{-1}))$ be the field of formal Laurent
  series over $\FF_q$, let $f\in\wh K$ be a quadratic irrational over
  $\FF_q(Y)$ and let $P\in R$ be an irreducible polynomial. We study
  the asymptotic properties of the degrees of the coefficients of the
  continued fraction expansion of quadratic irrationals such as $P^nf$
  as $n\ra+\infty$, proving, in sharp contrast with the case of
  quadratic irrationals in $\RR$ over $\QQ$ considered in
  \cite{AkaSha15}, that they have one such degree very large with
  respect to the other ones. We use arguments of \cite{BroPau07JLMS}
  giving a relationship with the discrete geodesic flow on the
  Bruhat-Tits building of $(\PGL_2,\wh K)$ and, with $A$ the diagonal
  subgroup of $\PGL_2(\wh K)$, the escape of mass phenomena of
  \cite{KemPauSha17} for $A$-invariant probability measures on the
  compact $A$-orbits along Hecke rays in the moduli space
  $\PGL_2(R)\bs\PGL_2(\wh K)$.
\end{abstract}

\section{Introduction}

This paper gives a positive characteristic analogue (with a
surprisingly different outcome) of results of Aka-Shapira
\cite{AkaSha15}. They studied the statistics of the period of the
continued fraction expansion of arithmetically constructed sequences
of quadratic irrational real numbers. For instance, given $\alpha\in
\RR$ a quadratic irrational over $\QQ$ and $p\in\ZZ$ a prime, they
proved that the equiprobability measures on the periodic part of the
orbit of the fractional part of $p^n\alpha$ under the Gauss map
weak-star converges towards the Gauss-measure on $[0,1]$ as
$n\ra+\infty$. Using the dynamical approach of \cite{KemPauSha17} and
giving the arithmetic applications announced there, we will prove that
in positive characteristic on the contrary, the coefficients of the
continued fraction expansion along similar sequences of quadratic
irrationals have a very irregular distribution (see in particular
Theorem \ref{theo:intro:distribution}), and that the proportion of
high degree coefficients might be positive.

Let us give a more precise description of our results. Let $\FF_q$ be
a finite field of order a positive power $q$ of a prime $p$ different
from $2$, and let $K=\FF_q(Y)$ be the field of rational functions in
one variable $Y$ over $\FF_q$. Let $R=\FF_q[Y]$ be the ring of
polynomials in $Y$ over $\FF_q$, let $\wh{K}=\FF_q((Y^{-1}))$ be the
non-Archimedean local field of formal Laurent series in $Y^{-1}$ over
$\FF_q$, let $\OOOO=\FF_q[[Y^{-1}]]$ be the local ring of $\wh{K}$
(consisting of formal power series in $Y^{-1}$ over $\FF_q$), and let
$\M=Y^{-1}\OOOO$ be the maximal ideal of $\OOOO$.

Any element $f\in \wh{K}$ may be uniquely written $f=[f]+\{f\}$ with
$[f]\in R$ (called the {\it integral part} of $f$) and $\{f\}\in \M$
(called the {\it fractional part} of $f$).  The {\it Artin map} $\Psi:
\M-\{0\}\ra \M$ is defined by $f\mapsto \big\{\frac{1}{f}\big\}$.  Any
$f\in \wh{K}$ which is irrational (that is, not in $K$) has a unique
continued fraction expansion
$$
f=a_0 +
\cfrac{1}{a_1+\cfrac{1}{a_2+ \cfrac{1}{a_3+\cdots}}}\;,
$$ 
with $a_0=a_0(f)=[f]\in R$ and $a_n=a_n(f)=
\big[\frac{1}{\Psi^{n-1}(f-a_0)} \big]$ a non constant polynomial for
$n\geq 1$ (called the {\it coefficients} of the continued fraction
expansion of $f$).

Let $QI=\{f\in \wh{K}\;:\;[K(f):K]=2\}$ be the set of quadratic
irrationals over $K$ in $\wh{K}$.  Given $f\in \wh{K}-K$, it is well
known that $f\in QI$ if and only if the continued fraction expansion
of $f$ is eventually periodic. Note that the projective action
(denoted by $\cdot$) of $\PGL_2(R)$ on $\PP_1(\wh K)=\wh K\cup
\{\infty\}$ preserves $QI$, keeping the periodic part of the continued
fraction expansions unchanged (up to cyclic permutation and invertible
elements, see \cite{BerNak00}), and that for all $x\in K$ and $f\in
QI$, we have $xf\in QI$. We refer for instance to the surveys
\cite{Lasjaunias00,Schmidt00} for background on the above notions.

\medskip
For every $f\in QI$, we denote by $\ell_f$ the length of the period of
the eventually periodic continued fraction expansion of $f$, and by
$d_{f,\,1},\dots, d_{f,\,\ell_f}\in\NN-\{0\}$ the degrees of the
polynomial coefficients appearing in this period.
For $c\in\;]0,1]$, we say that a sequence
$(f_n)_{n\in\NN}$ in $QI$ has {\it $c$-degree-escaping continued
fraction expansions} if for every $N\in\NN$,
$$
\liminf_{n\ra+\infty}\;\;\frac{\max_{i=1,\,\dots,\,\ell_{f_n}}d_{f_n,\,i} -N}
{\sum_{i=1}^{\ell_{f_n}} d_{f_n,\,i}}\geq c\;.
$$

\btheo\label{theo:escapingcfe}  Let $f\in QI$.
\begin{enumerate}
\item[(1)] 
  For every irreducible polynomial $P\in R$, there exists
  $c=c_{f,\,P}>0$ such that the sequence $(P^nf)_{n\in\NN}$ has
  $c$-degree-escaping continued fraction expansions.
\item[(2)] 
  Let $\ga\in \GL_2(R)$ projectively fixing $f$ whose
  discriminant $P=(\tr\ga)^2-4\det\ga$ is irreducible over
  $\FF_q$. Then for all $N\in\NN$,
\begin{equation}\label{eq:fullescape}
\limsup_{n\ra+\infty}\;\;
\frac{\max_{i=1,\,\dots,\,\ell_{P^nf}}d_{P^nf,\,i} -N}
{\sum_{i=1}^{\ell_{P^nf}} d_{P^nf,\,i}}=1\;.
\end{equation}
\item[(3)] 
  For every irreducible polynomial $P\in R$, there exists
  $c=c_{f,\,P}>0$ such that for uncountably many sequences
  $(\ga'_n)_{n\in\NN}$ in $\PGL(\FF_q[Y,\frac{1}{P}])$, the sequence
  $(\ga'_n\cdot f)_{n\in\NN}$ has $c$-degree-escaping continued fraction
  expansions.
\end{enumerate}
\etheo

Given $f\in QI$ and an irreducible\footnote{This irreducibility
  assumption in not essential, it only simplifies the exposition so
  that we don't have to work with several completions simultaneously.}
polynomial $P\in R$, the study of the properties of the continued
fraction expansion of $(P^nf)_{n\in\NN}$ is interesting in particular
by the analogy with the case of quadratic irrationals in $\RR$ over
$\QQ$, and the above results are strikingly different from the ones
obtained in \cite{AkaSha15}. Indeed, given $f\in\RR$ a quadratic
irrationnal over $\QQ$ whose continued fraction expansion
$(a_n(f))_{n\in\NN}$ is periodic after the time $k_f\geq 1$ with
period $\ell_f$, as a consequence of the much stronger results in
\cite{AkaSha15}, we have, for every positive prime integer $p$,
$$
\lim_{n\ra+\infty}\;\;
\frac{\max_{i=k_f,\,\dots,\,k_f+\ell_{p^nf}-1} \,\log \,a_{i}(p^nf)}
{\sum_{i=k_f}^{k_f+\ell_{p^nf}-1} \,\log \,a_{i}(p^nf)}\;=\;0\;.
$$

Theorem 4.5 of \cite{deMTeu04} has, using \cite[Lem.~4.2]{deMTeu04},
the following consequence : If $(a_{n,k})_{k\in\NN}$ is the continued
fraction expansion of $P_nf$ for all $n\in\NN$, where
$(P_n)_{n\in\NN}$ is a sequence in $\FF_q[Y]$ with $P_n$ dividing
$P_{n+1}$ and $P_{n+1}/P_n$ having bounded degree, then
$\sup_{n,k\in\NN-\{0\}} \deg a_{n,k}=+\infty$. Assertion (1) of
Theorem \ref{theo:escapingcfe} implies a much stronger result for the
sequence $(P_n=P^n)_{n\in\NN}$ of powers of any irreducible polynomial
$P\in R$.

Assertion (2) of Theorem \ref{theo:escapingcfe}, with the explicit
example given after \cite[Lem.~13]{KemPauSha17}, implies for instance
that if $p=q=3$, $P=Y^2+4$ and $f=\frac{Y-\sqrt{Y^2+4}}{2}$ so
that the Galois conjugate $f^\sigma$ of $f$ is the well-known
continued fraction expansion $f^\sigma=Y +
\cfrac{1}{Y+\cfrac{1}{Y+\cdots}}$, then Equation \eqref{eq:fullescape}
holds.

\medskip
The main point in the proof of Theorem \ref{theo:escapingcfe} is to
relate the orbits of quadratic irrationals under the Artin map to the
compact orbits of the diagonal subgroup $A$ of $G=\PGL_2(\wh K)$ on
the moduli space $X$ of $R$-lattices (up to homotheties) in $\wh
K\times \wh K$.  More precisely (forgetting about measure zero subsets
in this introduction, see Section \ref{sec:relatdiagArtin} for exact
statements), let $\Ga=\PGL_2(R)$ be Nagao's lattice (see
\cite{Nagao59,Weil70}) so that $X$ is the homogeneous space $\Ga\bs G$
endowed with the action of $A$ by translations on the right. We
construct, using \cite{BroPau07JLMS} and \cite[\S 9.6]{EinWar11} as
inspiration, a natural cross-section $C$ for the action of $A$ on $X$
and a natural map $\Theta_2: C\ra \M$ such that if $T:C\ra C$ is the
first return map in $C$ of the orbits under $A$, then we have a
commutative diagram
$$
\begin{array}{ccc}
C &\stackrel{T}{\longrightarrow}& C \medskip\\
^{\Theta_2}\downarrow\;\;\;\;& &\;\;\;\;\downarrow{} ^{\Theta_2}\medskip\\
\M &\stackrel{\Psi}{\longrightarrow}& \M\;.
\end{array}
$$ 
Given a quadratic irrational $f$, assuming for simplicity that
$f\in \M$ and that its Galois conjugate does not belong to $\OOOO$, we
construct in Lemma \ref{lem:xfAper} an explicit element $x_f\in C$
whose $A$-orbit in $X$ is compact, such that the image by $\Theta_2$
of the intersection with $C$ of the orbit of $x_f$ by right
translations under the semigroup $A_+=\Big\{ \begin{bmatrix} 1 & 0
  \\ 0 & f \end{bmatrix}\;:\;f\in\wh K-\OOOO\Big\}$, is equal to the orbit
of $f$ under the iterations of $\Psi$.

We then apply in Section \ref{sec:applyKPS} the results of
\cite{KemPauSha17} on the escape of mass of the $A$-invariant
probability measures on the compact $A$-orbits in $X$ varying along
rays in the Hecke trees of elements of $X$ associated to a
given irreducible polynomial in $R$. 

\medskip 
We conclude this introduction by stating a distribution result for the
orbits under the Artin maps of arithmeticaly defined families of
quadratic irrationals, again in sharp contrast with the ones obtained
in \cite{AkaSha15}, and again using the above-mentionned results of
\cite{KemPauSha17}. For every $f\in QI$, let $\nu_f$ be the
equiprobability on the periodic part of the $\Psi$-orbit of $\{f\}$.
Given nice sequences $(f_n)_{n\in\NN}$ in $QI$, it is interesting to
study the weak-star accumulation points of the equiprobabilities
$\nu_{f_n}$ as $n\ra+\infty$.

\btheo \label{theo:intro:distribution}
Let $f\in QI$ and $P\in R$ be an irreducible polynomial. For
uncountably many sequences $(\ga''_n)_{n\in\NN}$ in $\PGL_2(R)$, for
every $f'\in QI$, there exists $c>0$ and a weak-star accumulation
point $\theta$ of the equiprobabilities $\nu_{P^n\ga''_n\cdot f}$ such
that $\theta\geq c\, \nu_{f'}$.  
\etheo

In particular, $\theta$ is not absolutely continuous with respect to
the Haar measure of $\M$, which is the analog for the Artin map $\Psi$
of the Gauss measure for the Gauss map. Thus this proves that there
are many Hecke-type sequences $(P^{n_k}\ga_{n_k}\cdot f)_{k\in\NN}$ of
quadratic irrationals whose continued fraction expansions not only
have coefficients with degree almost equal to the sum of the degrees
on the period, but also have a positive proportion of coefficients
which are constant (say equal to $Y$, taking for instance $f' = Y+
\cfrac{1}{Y+\cfrac{1}{Y+\cdots}}$).

\bigskip\noindent{\small {\em Acknowledgements: } A first version of
  this work was partially supported by the NSF grants DMS-1440140,
  while the first author was in residence at the MSRI, Berkeley CA,
  during the Fall 2016 semester.  The first author would like to thank
  the Isaac Newton Institute for Mathematical Sciences, Cambridge, for
  support and hospitality during the programme ``Non-positive
  curvature group actions and cohomology'' where work on this paper
  was completed. This work was supported by EPSRC grant no
  EP/K032208/1. The second author is supported by ISF grants 357/13
  and 871/17.}

\section{A cross-section  for the orbits of the diagonal group}
\label{sec:relatdiagArtin}

In this first section, we give an explicit correspondence between the
orbits in the moduli space $\PGL_2(\Ga) \bs\PGL_2(\wh K)$ under the
diagonal subgroup and the orbits in $\M$ under the Artin map $\Psi$.
We start by recalling some notation and definitions from
\cite{Serre83, Paulin02, BroPau07JLMS, KemPauSha17}.

We denote by $v_\infty:K\ra\ZZ\cup\{+\infty\}$ the (normalized,
discrete) valuation on $K$ such that $v_\infty(P)= - \deg P$ for all
$P\in R-\{0\}$ (and $v_\infty(0)=+\infty$), and again by $v_\infty:\wh
K\ra\ZZ\cup\{+\infty\}$ its continuous extension to $\wh K$ for the
absolute value
$$
|x|_\infty = q^{-v_\infty(x)}\;.
$$ 
The convergence of the continued fraction expansions takes place in
$\wh K$ endowed with $|\cdot|_\infty$.

Let $G$ be the locally compact group $\PGL_2(\wh K)$, and let $\Ga=
\PGL_2(R)$, which is a nonuniform lattice in $G$. Let $A$ be the image
in $G$ of the diagonal subgroup of $\GL_2(\wh K)$.  We denote by
$\uA(\OOOO)$ the image in $G$ of the group of diagonal matrices with
diagonal coefficients in $\OOOO^\times$. We denote by
$\begin{bmatrix} a & b \\ c & d \end{bmatrix}$ the image in $G$ of
$\begin{pmatrix} a & b \\ c & d \end{pmatrix}\in \GL_2(\wh K)$, and by
$\alpha_0: \wh{K}^\times\ra A$ the group isomorphism $t\mapsto
\begin{bmatrix} 1 & 0 \\ 0 & t \end{bmatrix}$. From now on, we consider 
$\aaa_0=\alpha_0(Y)$.

Let $X= \Ga\bs G$. The map $\Ga g\mapsto g^{-1}[R\times R]$ identifies
$X$ with the space of homothety classes of $R$-lattices in
$\wh{K}\times \wh{K}$ (that is, of rank $2$ free $R$-submodules
spanning the vector plane $\wh{K}\times \wh{K}$ over $\wh{K}$).  A
point $x\in X$ is called {\it $A$-periodic} if its orbit by right
translations under $A$ is compact. This orbit $x A$ then carries a
unique $A$-invariant probability measure, denoted by $\mu_x$.

We denote by $\PP_1(\wh{K})^{(3)}$ the space of triples of pairwise
distinct points in $\PP_1(\wh{K})$, and by $\PP_1(\wh{K})^{(3)
  \natural}$ its intersection with $(\wh{K}-K)^2\times \PP_1(\wh{K})$.
Recall that the diagonal action of $G=\PGL_2(\wh K)$ on
$\PP_1(\wh{K})^{(3)}$ is simply transitive, so that $G$ identifies
with $\PP_1(\wh{K})^{(3)}$ by the map $g\mapsto g(\infty,0,1)$ and $X$
identifies with $\Ga\bs\PP_1(\wh{K})^{(3)} $ by the map $\Ga g\mapsto
\Ga g(\infty,0,1)$.  We denote by $G^\natural$ the dense measurable
subset of $G$ corresponding to $\PP_1(\wh{K})^{(3) \natural}$. It has
full measure for any Haar measure on $G$. It is $\Ga$-invariant on the
left, since the projective action of $\Ga$ preserves $\PP_1(\wh{K})
-\PP_1(K)=\wh{K}-K$. It is $A$-invariant on the right, since the
projective action of $A$ preserves $\infty$ and $0$. We denote by
$X^\natural$ the dense measurable subset of $X$ corresponding to
$\Ga\bs\PP_1(\wh{K})^{(3) \natural}$, which is $A$-invariant on the
right and has full measure for the $G$-invariant probability measure
$m_X$ on $X$. We also define
$$
\M^\natural=\M-(\M\cap K) \;\;\;{\rm and}\;\;\; 
(^c \OOOO)^\natural=\wh{K}-(\OOOO\cup K)\;,
$$ 
the subsets of irrational points in $\M$ and $^c\OOOO$ respectively.

\medskip
Let us recall the minimal necessary background on the {\it Bruhat-Tits
  tree} $\TT$ of $(\PGL_2, \wh{K})$ from \cite{Serre83} or \cite[\S
  2.3]{KemPauSha17}.  It is a regular tree of degree $q+1$, endowed
with an action of $G$ by automorphisms (without inversions). Its
(discrete) set of vertices $V\TT$ is the set of homothety classes
(under $\wh{K}^\times$) of $\OOOO$-lattices in $\wh{K}\times \wh{K}$,
endowed with the maximal distance such that the distance between the
endpoints of any edge is $1$.

We denote by $E\TT$ its set of edges, by $o:E\TT\ra V\TT$ the origin
map of edges, and by $*=[\OOOO\times \OOOO]$ its standard base point,
whose stabilizer in $G$ is the maximal compact subgroup
$\PGL_2(\OOOO)$ of $G$, so that the map $g\mapsto g*$ induces an
identification between $G/\PGL_2(\OOOO)$ and $V\TT$.

The boundary at infinity $\partial_\infty \TT$ of $\TT$ identifies
($G$-equivariantly) with the projective line $\PP_1(\wh{K})$ over
$\wh{K}$, which itself identifies with $\wh{K} \cup \{\infty\}$ using
the map $[x:y] \mapsto xy^{-1}$.  We endow it with the unique
probability measure $\mu_{\rm Hau}$ invariant under the compact group
$\PGL_2(\OOOO)$. We denote by $g\cdot\xi$ the projective action of
$g\in G$ on $\xi\in \PP_1(\wh{K})=\partial_\infty \TT$.

We say that an edge $e$ of $\TT$ {\it points towards} a point at
infinity $\xi\in\partial_\infty \TT$ if its terminus $t(e)$ belongs to
the geodesic ray starting from its origin $o(e)$, with point at
infinity $\xi$.

We denote by $\G\TT$ the space of {\it geodesic lines} in $\TT$ (that
is, the set of isometric maps $\ell:\ZZ\ra V\TT$ endowed with the
compact-open topology), and by $\ell_*\in\G\TT$ the unique geodesic
line with $\ell_*(-\infty) =\infty$, $\ell_*(+\infty) = 0$, and
$\ell_*(0)=*$, so that $\ell_*(n)=\aaa_0^{\,n}*$ for every
$n\in\ZZ$. We denote by $(\phi_n)_{n\in\ZZ}$ the (discrete time) {\it
  geodesic flow} on $\G\TT$, where $\phi_n\ell:k\mapsto \ell(k+n)$ for
all $n\in\ZZ$ and $\ell \in \G\TT$, as well as its quotient flow on
$\Ga\bs\G\TT$. The stabilizer of $\ell_*$ for the transitive action
$(g,\ell)\mapsto \{g\,\ell:k\mapsto g\,\ell(k)\}$ of $G$ on $\G\TT$ is
exactly $\uA(\OOOO)$. Hence the map $g\mapsto g\,\ell_*$ induces a
homeomorphism from $G/\uA(\OOOO)$ to $\G\TT$. We have the following
crucial property relating the right action of $A$ on $X=\Ga\bs G$ and
the geodesic flow on $\G\TT$: For all $g\in G$ and $n\in\ZZ$,
\begin{equation}\label{eq:commuta0geodflow}
\phi_n (g\,\ell_*)= g\,\aaa_0^{\,n}\,\ell_*\;.
\end{equation}

\medskip
Recall that the {\it Buseman function} of $\TT$ is the map
$\beta:\partial_\infty\TT\times V\TT\times V\TT\ra\RR$ defined by
$$ 
\beta:(\xi,x,y)\mapsto \beta_\xi(x,y)=d(x,p)-d(y,p)
$$ 
for any point $p$ close enough to $\xi$, and in particular the point
$p\in V\TT$ such that the geodesic rays $[x,\xi[$ and $[y,\xi[$ with
origin respectively $x$ and $y$, with point at infinity $\xi$, intersect
exactly in the geodesic ray $[p,\xi[\,$. A (closed) {\it horoball} in
$\TT$ {\it centered at} a point at infinity $\xi\in\partial_\infty
\TT$ is a sublevel set $\{x\in V\TT\;:\;\beta_\xi(x,x_0)\leq 0\}$ for
some $x_0\in V\TT$.

We denote (see for instance \cite[\S 6.2]{Paulin02}) by
$(\H_\xi)_{\xi\in\PP_1(K)}$ the $\Ga$-invariant family of maximal
horoballs in $\TT$ with pairwise disjoint interiors such that
the point at infinity of $\H_\xi$ is $\xi$. If $\H_\infty = \{x\in
V\TT\;:\;\beta_\infty(x,*)\leq 0\}$ is the horoball centered at
$\infty$ whose boundary contains $*$, then $\H_\xi=\ga\H_\infty$ for
every $\ga\in\Ga$ such that $\ga\cdot\infty=\xi$. Note that
$$
\Ga\, * =\bigcup_{\xi\in\PP_1(K)}\partial\H_\xi= 
V\TT-\bigcup_{\xi\in\PP_1(K)}\stackrel{\circ}{\H_\xi}\;.
$$ 
For instance (see loc.cit.), the horoballs in this family whose
boundary contains the base point $*$ are the one centered at the
elements in $\PP_1(\FF_q)=\FF_q\cup\{\infty\}$, and $\M$ (respectively
$^c\OOOO$, and $\OOOO^\times$) is the set of points at infinity of the
geodesic rays starting from $*$ whose first edge lies in $\H_0$
(respectively in $\H_\infty$, and neither in $\H_0$ nor in
$\H_\infty$).

\medskip
Let 
\begin{align*}
\wt C&=\{g\in G^\natural\;:\;\; 
g\cdot \infty\in \;^c\OOOO,\;\;g\cdot 0\in\M,\;\;g\cdot 1=1\}\\ 
&=\{g\in G\;:\;\; g\cdot \infty\in \;(^c \OOOO)^\natural,
\;\;g\cdot 0\in\M^\natural,\;\;g\cdot 1=1\}\;.
\end{align*} 
As we shall see in the proof below, we actually have $\wt C\subset
\PGL_2(\OOOO)$ and if $g\mapsto \ov{g}$ is the group morphism
$\PGL_2(\OOOO)\ra \PGL_2(\OOOO/\M)=\PGL_2(\FF_q)$ of reduction modulo
$\M$, then $\ov{g}=\id$ for every $g\in \wt C$.

Let $C=\Ga\wt C$ be the image of $\wt C$ in $X^\natural=\Ga\bs
G^\natural$ by the canonical projection. We will in particular prove
in the next result that $C$ is a natural cross-section for the right
action of $A$~: Every orbit of $A$ in $X^\natural$ meets $C$ in an
infinite countable subset.\footnote{This follows from Assertion (3)
  and Assertion (4) below.} As in \cite[Theo.~3.7]{BroPau07JLMS}, let
$$ 
\wt\Psi:\;({}^c\OOOO)^\natural\;\times\; \M^\natural\ra
\;({}^c\OOOO)^\natural\;\times\; \M^\natural
$$ 
be the map defined by
$$
(\xi_-,\xi_+)\mapsto
\Big(\frac{1}{\xi_-}-\Big[\frac{1}{\xi_+}\Big],\Psi(\xi_+)\Big)\;.
$$

\btheo\label{theo:commutfrmArtin} \begin{enumerate}
\item[(1)] The subset $C$ of $X^\natural$ is closed in $X^\natural$,
  and the restriction to $\wt C$ of the canonical projection $G\ra
  X=\Ga\bs G$ is a homeomorphism from $\wt C$ to $C$.
\item[(2)] The map $\wt \Theta: g\mapsto (g\cdot \infty, g\cdot 0)$ from
  $\wt C$ to $\PP_1(\wh K)^2$ induces a homeomorphism
$$
\Theta: C\ra ({}^c\OOOO)^\natural\;\times\; \M^\natural\;.
$$
\item[(3)] Every $A$-orbit in $X^\natural$ meets $C$. 
\item[(4)] For every $x\in
  C$, there exist unique sequences $(\alpha_k=\alpha_k(x))_{k\in\ZZ}$ in
  $\uA(\OOOO)$ and $(t_k=t_k(x))_{k\in\ZZ}$ in $2\ZZ$ with $\alpha_0=\id$,
  $t_0=0$ and $t_k<t_{k+1}$ for all $k\in\ZZ$ such that
$$
xA\cap C=\{x\alpha_k{\aaa_0}^{t_k}\;:\;k\in\ZZ\}\;.
$$
We denote by $T:C\ra C$ the map $x\mapsto x\alpha_1{\aaa_0}^{t_1}$.
\item[(5)] The following diagram is commutative
$$
\begin{array}{ccc}
 C &\stackrel{\Theta}{\longrightarrow}& 
({}^c\OOOO)^\natural\;\times\; \M^\natural\;{}\medskip\\
T\downarrow\;\;\;& &\;\downarrow \wt\Psi\medskip\\
 C &\stackrel{\Theta}{\longrightarrow}& 
({}^c\OOOO)^\natural\;\times\; \M^\natural\;.\\
\end{array}
$$
\end{enumerate}
\etheo

We call $C$ the {\it cross-section} for the action of $A$ on
$X^\natural$. We say that $t_1(x)$ is the {\it first return time} in
$C$ of the orbit of $x\in C$ under $A$, and we call $T$ the {\it first
  return map} of the orbits under $A$ in the cross-section $C$.

Note that $T$ is a homeomorphism, with $t_m(T^nx)=t_{m+n}(x)-t_n(x)$ and 
$\alpha_m(T^nx)=\alpha_{m+n}(x)\alpha_n(x)^{-1}$ for all $m,n\in\ZZ$.

The commutativity of the diagram in Assertion (5) is inspired by, and
closely related to, the commutativity of the diagram in
\cite[Theo.~3.7]{BroPau07JLMS}.

\medskip
\dem (1) The closure $\overline{\wt C}$ of $\wt C$ in $G$ is $\{g\in
G\;:\; g\cdot\infty\in\;^c\OOOO, g\cdot 0\in\M , g\cdot 1=1\}$, which
is closed since the conditions on $g\cdot\infty$, $g\cdot 0$ and
$g\cdot 1$ are closed. By the description of the geodesic rays in
$\TT$ starting from $*$, any geodesic line $\ell$ with endpoints
$\ell(-\infty)\in \;^c\OOOO$ and $\ell(+\infty)\in\M$ passes through
the base point $*$, and the closest point projection of $1\in\PP_1(\wh
K)= \partial_\infty\TT$ on $\ell$ is equal to $*$. Hence if $g\in
\overline{\wt C}$, since $g$ commutes with the closest point
projections on geodesic lines, we have $g*=*$. Therefore $g\in
\overline{\wt C}$ is contained in the stabilizer $\PGL_2(\OOOO)$ of
$*$ in $G$, which is compact.  Hence $\overline{\wt C}$ is compact and
so is its image $\Ga \overline{\wt C}$ in $X=\Ga\bs G$. Thus
$C=X^\natural \cap \Ga \overline{\wt C}$ is closed in $X^\sharp$.

Let us prove that the canonical projection $G\ra X$ is injective on
$\overline{\wt C}$. Since its image in $X$ is $\Ga \overline{\wt C}$,
and since $\overline{\wt C}$ is compact and $X$ is Hausdorff, this
proves that the orbit map $\overline{\wt C}\ra\Ga \overline{\wt C}$ is
a homeomorphism, and so is $\wt C\ra C$ by restriction.

Let us denote by $g\mapsto \ov{g}$ the group morphism $\PGL_2(\OOOO)
\ra \PGL_2(\OOOO/\M)=\PGL_2(\FF_q)$ of reduction of coefficients
modulo $\M$, which is the identity map on $\PGL_2(\FF_q)$. Recall (see
for instance \cite{Serre83}) that we have an identification between
the {\it link} 
$$
lk(*)=\{e\in E\TT\;:\;o(e)=*\}
$$ 
of the vertex $*$ in $\TT$ and the finite projective line
$\PP_1(\FF_q) = \FF_q\cup \{\infty\}$ which is equivariant under the
reduction morphism $\PGL_2(\OOOO)\ra \PGL_2(\FF_q)$. For every
$\ga\in\Ga$, if there exists $g,g'\in \overline{\wt C}$ such that $\ga
g=g'$, then since $\OOOO\cap R=\FF_q$, we have 
$$
\ga= g'g^{-1}\in\PGL_2(\OOOO)\cap \Ga=\PGL_2(\FF_q)\;.
$$ 
Furthermore, $\ga$ fixes the edges with origin $*$ pointing towards
$\infty$, $0$ and $1$ respectively, hence $\ga$ is an element of
$\PGL_2(\FF_q)$ fixing the points $\infty$, $0$ and $1$ in
$\PP_1(\FF_q)$. Since $\PGL_2(\FF_q)$ acts simply transitively on the
set of triples of pairwise distinct elements of $\PP_1(\FF_q)$, we
have $\ga=\id$. Therefore the canonical projection $G\ra \Ga\bs G$ is
indeed injective on $\overline{\wt C}$. This completes the proof of
Assertion (1).

\medskip \noindent (2) The map from $G$ to $\PP_1(\wh K)^{(3)}$
defined by $g\mapsto (g\cdot\infty,g\cdot 0,g\cdot 1)$ being a
homeomorphism, the map from $\overline{\wt C}$ to $\M\times \;^c\OOOO$
defined by $g\mapsto (g\cdot\infty,g\cdot 0)$ is a
homeomorphism. Hence by restriction, the map from $\wt C$ to
$\M^\natural\times (^c\OOOO)^\natural$ defined by $g\mapsto
(g\cdot\infty,g\cdot 0)$ is a homeomorphism. Assertion (2) then
follows from Assertion (1), by composition.

\medskip \noindent (3) Since $A$ is abelian, the direct product group
$\Ga\times A$ acts on $G$ by $((\ga,a),g)\mapsto \ga g a$ for all
$\ga\in\Ga,g\in G,a\in A$, and this action preserves
$G^\natural$. Let us prove that every $(\Ga\times A)$-orbit in
$G^\natural$ meets $\wt C$. By taking quotients on the left modulo
$\Ga$, this proves that every $A$-orbit in $X^\natural$ meets $C$.

Let us fix $g\in G^\natural$. Let us prove that $\Ga g A$ meets $\wt
C$. Since the family of horoballs $(\H_\xi)_{\xi\in\PP_1(K)}$
covers $\TT$, since its elements have pairwise disjoint interiors,
since the boundary of $\H_\xi$ is contained in $\Ga *$ for every
$\xi\in\PP_1(K)$, and since no complete geodesic line is contained in
a horoball, the geodesic line $g\ell_*$ meets $\Ga *$. Hence, using
Equation \eqref{eq:commuta0geodflow}, there exists $n\in\ZZ$ such that
$g \aaa_0^{\, n} * = g\ell_*(n)\in\Ga*$. We may hence assume, up to
multiplying $g$ on the right by an element of $\aaa_0^\ZZ\subset A$, that
$g\,*\in \Ga *$.

Up to multiplying $g$ on the left by an element of $\Ga$, we may hence
assume that $g *=*$. Since the link $lk(*)$ of $*$ identifies with the
projective line $\PP_1(\FF_q)$, since $\PGL_2(\FF_q)$ is contained in
$\Ga$ and acts transitively on the ordered pairs of distinct elements
of $\PP_1(\FF_q)$, we may assume, up to multiplying $g$ on the left by
an element of $\Ga$, that the edge $e_\infty$ with origin $*$ pointing
towards $\infty$ also points towards $g\cdot\infty$ and that the edge
$e_0$ with origin $*$ pointing towards $0$ also points towards $g\cdot
0$.

By for instance \cite{Paulin02} and as recalled before the definition
of $\wt C$, a geodesic ray starting from $*$ has its point at infinity
in respectively $^c\OOOO$ or $\M$ if and only if its first edge is
$e_\infty$ or $e_0$. In particular, $g\cdot
\infty\in({}^c\OOOO)^\natural$ and $g\cdot 0\in\M^\natural$.

\begin{center}
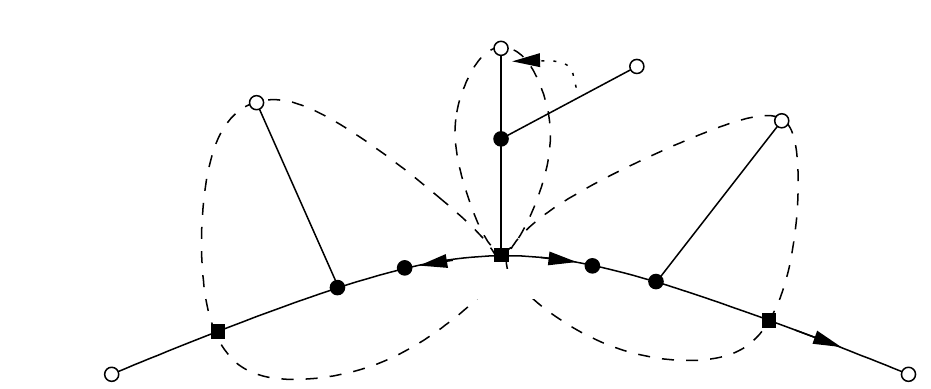
\end{center}

If the projection of a point at infinity to the connected union of two
edges is their middle vertex, then the projection of a point at
infinity to any geodesic line containing these two edges is again this
middle vertex. Note that projections on closed subtrees commute with
tree automorphisms. Hence the projection of $g\cdot 1$ on the geodesic
line $\ell_*$ with points at infinity $\infty$ and $0$ (which contains
$e_0$ and $e_\infty$) is equal to $*$.  The group $\uA(\OOOO)$ acts
simply transitively on $\OOOO^\times$, which is the subset of
$\partial_\infty\TT=\wh K\cup \{\infty\}$ consisting of elements whose
projections on $\ell_*$ is $*$. Hence there exists $\beta_0\in
\uA(\OOOO)$ such that $g\beta_0\cdot 1=1$. Since multiplying $g$ on
the right by $\beta_0$ does not change $g\cdot \infty$ and $g\cdot
0$, we hence have $g\beta_0\in \wt C$. This proves Assertion (3).

\medskip \noindent (4) Let us now fix $x\in C$, and let $g\in\wt C$ be
the unique element in $\wt C$ such that $x=\Ga g$.

Recall that any geodesic line which enters into the interior of a
(closed) horoball and does not converge to its point at infinity has
to exit this horoball after a finite time, and that the distance
travelled inside the horoball is even. Since $g\cdot \infty$ and
$g\cdot 0$ do not belong to $K$, there exists a unique sequence
$(t_k=t_k(x))_{k\in\ZZ}$ in $2\ZZ$ such that $t_0=0$, $t_k<t_{k+1}$
for all $k\in\ZZ$ and
$$
g\ell_*(\ZZ)\cap (\Ga *)=\{g\ell_*(t_k)\;:\;k\in\ZZ\}\;.
$$ 
This sequence is hence the sequence of times at which the geodesic
line $g\ell_*$ passes through $\Ga *$ (normalized by being at time
$t_0=0$ at $*$).

Let us fix $a\in A$ such that $x\,a\in C$, and let us prove that there
exists a unique pair $(k,\alpha)\in\ZZ\times \uA(\OOOO)$ such that $a=
\alpha \,\aaa_0^{t_k}$.

\begin{center}
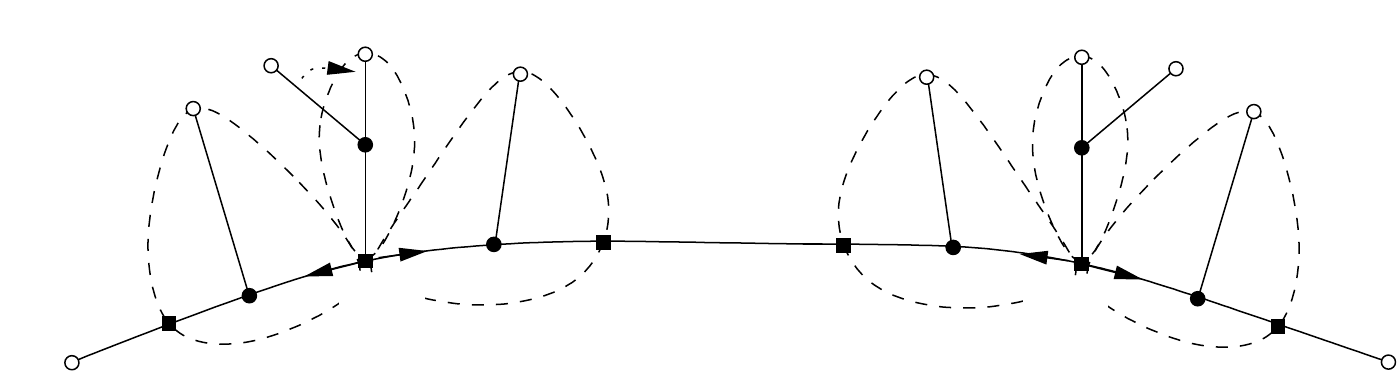
\end{center}

We have $g\,a\,*\in\Ga\,*$ and $g\,a\cdot 1$ is a rational point at
infinity. Indeed, since the lift of $x\,a$ in $\wt C$ has the form
$\ga \,g\,a$ for some $\ga\in\Ga$, we have $g\,a\, *= \ga^{-1} \,*\in
\Ga\, *$ and $g\,a\cdot 1= \ga^{-1}\cdot 1\in K\cup\{\infty\}$. In
particular, there exists a unique $k\in\ZZ$ such that $g\,a\, *=
g\,\aaa_0^{t_k}\,*$.

Let $\xi_{k-1}$ and $\xi_k$ be the points at infinity of the
horoballs in the family $(\H_\xi)_{\xi\in\PP_1(K)}$ into which the
edges $g\,a\,e_\infty$ and $g\,a \,e_0$ respectively enter.

By the simple transitivity of the projective action of $\PGL_2(K)$ on
$\PP_1(K)^{(3)}$, let $\ga\in\Ga$ be the unique element of $\Ga$
sending $\xi_{k-1},\xi_k,g\,a\cdot 1$ to respectively $\infty, 0,1$.
Then $\ga$ sends $g\,a\,*$ to $*$, $g\,a\,e_\infty$ to $e_\infty$ and
$g\,a \,e_0$ to $e_0$. Therefore $\ga\,g\,a\cdot\infty$, which is
irrational since equal to $\ga\,g\cdot\infty$ and $g\cdot\infty$ is
irrational, is the point at infinity of a geodesic ray starting by the
edge $e_\infty$, hence belongs to $({}^c\OOOO)^\natural$. Similarly,
$\ga\,g\,a\cdot\infty$ belongs to $\M^\natural$. As $\ga\,g\,a\cdot
1=1$ by construction of $\ga$, we have $\ga\,g\,a\in \wt C$, and by
Assertion (1),  $\ga$ is the unique element of $\Ga$ such that
$\ga\,g\,a\in \wt C$.

Using arguments similar to the ones at the end of the proof of
Assertion (3), since $g\,\aaa_0^{t_k}\cdot 1$ projects to $g\,a\,*=
g\,\aaa_0^{t_k}*$ on the union of $g\,a\,e_\infty$ and $g\,a \,e_0$, the
point at infinity $\ga g\,\aaa_0^{t_k}\cdot 1$ projects to $\ga
\,g\,a\,*=*$ on the union of $\ga \,g\,a\,e_\infty=e_\infty$ and $\ga
\,g\,a \,e_0=e_0$. Hence there exists one and only one $\alpha\in
\uA(\OOOO)$ such that $\ga \,g\,\aaa_0^{t_k}\alpha\cdot 1=1$. Since $\ga
\,g\,\aaa_0^{t_k}\alpha$ and $\ga \,g\,a$ give the same images to
$\infty,0,1$, they are equal, hence $a=\aaa_0^{t_k}\alpha$, for a unique
$\alpha$. This proves Assertion (4).

\medskip 
\noindent 
(5) Let us fix $x\in C$. Let $n=t_1(x)$ be the first return time of
the $A$-orbit of $x$ to $C$, and $k=\frac{n}{2}$.  We denote by $g\in
\wt C$ the lift of $x$ to $\wt C$. In particular, $g$ fixes $*, 1, e_\infty,
e_0$.

\begin{center}
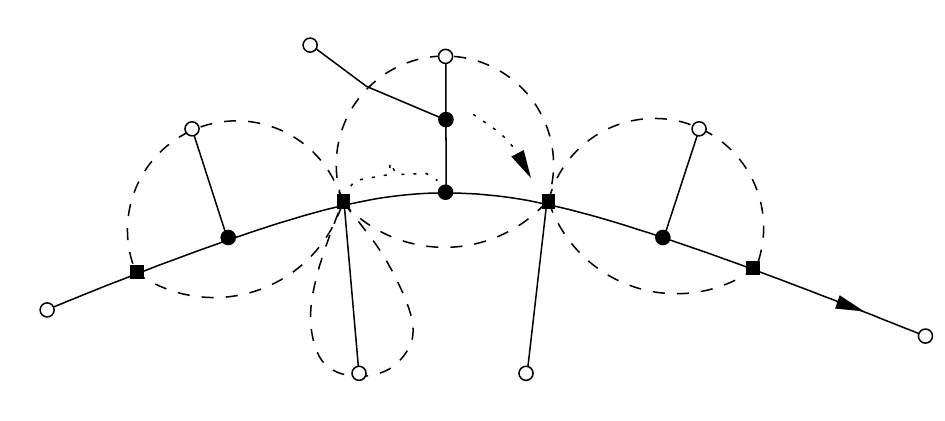
\end{center}

By the geometric interpretation of the continued fraction expansion of
$g\cdot 0\in \M$ (see \cite[\S 6]{Paulin02}), we have

$\bullet$~ 
the exiting time of $g\,\ell_*$ out of the horosphere $\H_0$ is $n$, and 
$$
k=\;\deg a_1(g\cdot 0)\;;
$$

$\bullet$~ the point at infinity of the horoball in the family
$(\H_\xi)_{\xi \in\PP_1(K)}$ into which enters the geodesic line
$g\,\ell_*$ after time $n$ is the first partial quotient
$$
\xi_1=\frac{1}{a_1(g\cdot 0)}
$$ 
of the continued fraction expansion of $g\cdot 0$. 

Note that $g \,\aaa_0^n$ maps $*$ to the unique intersection point of
$\H_0$ and $\H_{\xi_1}$, as well as $e_\infty$ and $e_0$ to the edges
starting from $g\, \aaa_0^n*$ and pointing towards the points at infinity
$0$ and $\xi_1$ respectively.

Define $\ga=\begin{bmatrix} -a_1(g\cdot 0) & 1 \\ 1 &
0 \end{bmatrix}$, which belongs to $\Ga$. Note that $\ga$, whose
associated homography is $z\mapsto \frac{1}{z}-a_1(g\cdot 0)$, maps
$0$ to $\infty$ and $\xi_1$ to $0$. Hence it maps the unique
intersection point of $\H_0$ and $\H_{\xi_1}$, which is
$g\,\aaa_0^n*$, to the unique intersection point of $\H_\infty$ and
$\H_{0}$, which is $*$. It also sends the first edge of the geodesic
ray from $g\,\aaa_0^n*$ to $0$, which is $g\, \aaa_0^ne_\infty$, to
the first edge of the geodesic ray from $*$ to $\infty$, which is
$e_\infty$. Similarly, $\ga \,g\,\aaa_0^n e_0=e_0$.

Since $\ga \,g\,\aaa_0^n$ fixes $e_\infty$ and $e_0$, the point at
infinity $\ga \,g\,\aaa_0^n\cdot \infty$ (which is equal to $\ga
\,g\cdot \infty$ hence is irrational) belongs to $^c\OOOO$, and
similarly $\ga \,g\,\aaa_0^n\cdot 0\in\M^\natural$. Furthermore, as
seen above, there exists a unique $\alpha \in\uA(\OOOO)$ such that
$\ga \,g\,\aaa_0^n\alpha \cdot 1=1$.

Therefore $\ga \,g\,\aaa_0^n\alpha\in\wt C$, by uniqueness $T(x)= x
\aaa_0^n\alpha$ and the unique lift of $T(x)$ in $\wt C$ is $\ga
\,g\,\aaa_0^n\alpha$.

We now compute
\begin{align*}
\Theta \circ T(x)&=\wt \Theta(\ga\,g\,\aaa_0^n\alpha)=
(\ga\,g\,\aaa_0^n\alpha\cdot\infty,\ga\,g\,\aaa_0^n\alpha\cdot 0)= 
(\ga\,g\cdot\infty,\ga\,g\cdot 0)\\ & =
\big(\frac{1}{g\cdot\infty}-a_1(g\cdot 0),\frac{1}{g\cdot 0}-
a_1(g\cdot 0)\big)= \wt \Psi (g\cdot\infty,g\cdot 0)=
\wt \Psi \circ \Theta(x)\;.
\end{align*}
This concludes the proof of Theorem \ref{theo:commutfrmArtin}.
\cqfd

\bigskip Let us introduce some notation before stating our next
result. All measures in this paper are Borel and nonegative.

We denote by $\mu_A$, $\mu_{\uA(\OOOO)}={\mu_A}_{\,\mid \uA(\OOOO)}$,
$\mu_{\wh K}$ and $\mu_{\M}={\mu_{\wh K}}_{\,\mid \M}$, the Haar
measures of the locally compact abelian groups $A$, $\uA(\OOOO)$,
$(\wh K,+)$ and $(\M,+)$ respectively, normalized as usual so that
$\mu_A(\uA(\OOOO))=1$ and $\mu_{\wh K}(\OOOO)=1$.

For every $x\in C$, we define $x_\pm$ by
$$
\Theta (x)=(x_-,x_+)\in ({}^c\OOOO)^\natural\times\M^\natural\;.
$$
Let $\Theta_2: C\ra \M^\natural$ be the map $\Theta_2=
\pr_2\circ\;\Theta$, where $\pr_2:(\xi_-,\xi_+)\mapsto \xi_+$ is the
projection on the second factor, which is the continuous and
surjective map on $C=\Ga\wt C$ induced by the map $g\mapsto g\cdot 0$
from $\wt C$ to $\M^\natural$. By Theorem \ref{theo:commutfrmArtin} (5),
the following diagram commutes:
\begin{equation}\label{commutTheta2}
\begin{array}{ccc}
 C &\stackrel{T}{\longrightarrow}& 
C\;{}\medskip\\
^{\Theta_2}\downarrow\;\;\;& &\;\downarrow {}^{\Theta_2}\medskip\\
 \M^\natural &\stackrel{\Psi}{\longrightarrow}& 
\M^\natural\;.\\
\end{array}
\end{equation}

Since $\uA(\OOOO)$ is compact and open in $G$ and by the uniqueness
properties in Theorem \ref{theo:commutfrmArtin} (4), the map from
$C\times \uA(\OOOO)$ to $X$ defined by $(x,a)\mapsto xa$ is a
homeomorphism onto its image, which is an open neighborhood of $C$ in
$X$. Denoting by $C\uA(\OOOO)$ its image, the map $C\uA(\OOOO)\ra C$
defined by $xa\mapsto x$ for all $x\in C$ and $a\in \uA(\OOOO)$ is a
fibration with fibers the right orbits of $\uA(\OOOO)$. The product
group $T^\ZZ\times\uA(\OOOO)$ acts in the obvious way on
$C\uA(\OOOO)$. We still denote by $t_k:C\uA(\OOOO)\ra\ZZ$ the
$k$-th return time in $C\uA(\OOOO)$ of the orbits under $\aaa_0^\ZZ$,
that is $t_k(xa)= t_k(x)$ for all $x\in C$ and $a\in\uA(\OOOO)$.
Given any measure $\mu$ on $C$, we denote by $\mu\otimes
\mu_{\uA(\OOOO)}$ the pushforwards by $(x,a)\mapsto xa$ of the product
measure of $\mu$ and $\mu_{\uA(\OOOO)}$, which is a finite nonzero
measure on $C\uA(\OOOO)$ if $\mu$ is finite and nonzero.

We denote by $M_A$ the space of finite nonzero $A$-invariant measures
on $X^\natural$, by $M_{T^\ZZ\uA(\OOOO)}$ the space of finite nonzero
$T^\ZZ\times\uA(\OOOO)$-invariant measures $m$ on $C\uA(\OOOO)$ such
that the integral $\int_{C\uA(\OOOO)}t_1\,dm$ is finite and nonzero,
by $M_{T}$ the space of $T$-invariant measures $m$ on $C$ such that
$\int_{C}t_1\,dm$ is finite and nonzero, and by $M_{\M}$ the space of
$\Psi$-invariant measures $m$ on $\M$ such that $\int_{f\in\M} \deg
(a_1(f))\,dm(f)$ is finite and nonzero. We endow these spaces with
their weak-star topologies.

We have maps
$$
M_A\;\stackrel{F_1}{\longrightarrow}\; M_{T^\ZZ\uA(\OOOO)}
\;\stackrel{F_2}{\longleftarrow}\;M_{T}
\;\stackrel{(\Theta_2)_*}{\longrightarrow}\;M_{\M}\;,
$$ 
where $F_1:m\mapsto m_{\mid C\uA(\OOOO)}$, $F_2:\mu\mapsto \mu\otimes
\mu_{\uA(\OOOO)}$, and we define (see the proof below for the
invertibility of $F_2$)
$$
F=(\Theta_2)_*\circ(F_2)^{-1}\circ F_1:M_A\ra M_{\M}\;.
$$

We denote by $\delta_x$ the unit Dirac mass at any point $x$ in any
measurable space.

\medskip
The following result strengthens the statement announced in the
Introduction relating $A$-orbits to $\Psi$-orbits, and gives the
measure-theoretic correspondence. The use of semi-groups comes from
the fact that $\Psi$ is non invertible. Recall that the analogue, for
the Artin map $\Psi:x\mapsto \frac{1}{x}-\big[\,\frac{1}{x}\,\big]$,
of the Gauss measure on $[0,1]$ for the Gauss map $x\mapsto
\frac{1}{x}-\big\lfloor\,\frac{1}{x}\,\big\rfloor$, is the Haar
probability measure of the additive group $\M$ (see for instance
\cite[Prop.~5.6]{BroPau07JLMS} for a dynamical explanation).

\btheo\label{theo:troispoint} 
\begin{enumerate}
\item[(1)] 
For every $x\in C$, the image by $\Theta_2$ of the intersection
with $C$ of the orbit of $x$ by right translations under the
semigroup $A_+=\uA(\OOOO)\,\aaa_0^{\NN}$ is equal to the $\Psi$-orbit
of $\Theta_2(x)\,$:
$$
\Theta_2\big(C\cap (x\,A_+)\big)\;=\;
\Psi^\NN(\Theta_2(x))\;.
$$ 
In particular, if $x\in C$ is $A$-periodic, then $\Theta_2(x)$ is a
quadratic irrational.  Furthermore, the period of $x\in C$ under the
iteration of  $T$ is equal to the period of the periodic part of the
continued fraction expansion of $\Theta_2(x)$.
\item[(2)] The map $F:M_A\ra M_{\M}$ is a homeomorphism for the
  weak-star topologies, such that $F(c\,m +m')=c\,F(m) + F(m')$ and
  $F(m)\leq F(m'')$ for all $c>0$ and $m,m', m''\in M_A$ such that
  $m\leq m''$.
\item[(3)] 
  The restriction to $C\uA(\OOOO)$ of the $G$-invariant
  probability measure $m_X$ on $X$ disintegrates by the fibration
  $C\uA(\OOOO)\ra C$, with conditional measures the
  $\uA(\OOOO)$-invariant probability measures on the fibers, over a
  finite measure $m_C$ on $C$ such that, for every $x\in C$,
  using the homeomorphism $\Theta: x\mapsto (x_-,x_+)$,
$$
dm_C(x)= \frac{q(q-1)^2}{2} \;
\frac{d\mu_{\wh K}(x_-)}{|x_-|_\infty^4}d\mu_{\wh K}(x_+)\;.
$$ 
Furthermore, the pushforwards under $\Theta_2$ of $m_C$ is equal to
a multiple of the Haar measure on $\M$:
$$
(\Theta_2)_*m_C= \frac{q^4(q-1)^2}{2(q^2+q+1)} \;\mu_{\M}\;,
$$ 
and the measure on $\M$ corresponding to the $G$-invariant
probability measure $m_X$ on $X$ is a multiple of the Haar measure on
$\M$:
$$ 
F(m_X) = \frac{q^4(q-1)^2}{2(q^2+q+1)} \;\mu_{\M}\;.
$$
\item[(4)] If $x\in C$ is $A$-periodic in $X$, then for $k$ big
  enough, the measure on $\M$ corresponding to the $A$-invariant
  probability measure $\mu_{x}$ on the $A$-orbit of $x$,
  normalized to be a probability measure, is equal to the
  equiprobability on the (finite) $\Psi$-orbit of $\Theta_2(T^kx)$:
  More precisely, if $n_x$ is the period of the periodic part of the
  continued fraction expansion of $\Theta_2(x)$, then
$$ 
F(\mu_{x}) = 
\frac{1}{2\sum_{n=k}^{k+n_x-1} \deg a_n(\Theta_2(x))}\;
\sum_{0\leq n\leq n_x-1} \delta_{\Psi^{n+k}(\Theta_2(x))}\;.
$$
\end{enumerate}
\etheo

In the last claim, we may take $k=0$ if $x$ is $A_+$-periodic, and in
general $k$ may be taken to be any time after which the continued
fraction expansion $\big(a_n(\Theta_2(x))\big)_{n\geq 1}$ of
$\Theta_2(x)$ is periodic.

\medskip
\dem (1) It follows from Theorem \ref{theo:commutfrmArtin} (4) that
$$
C\cap (x\,A_+)= \{x\alpha_k{\aaa_0}^{t_k}\;:\;k\in\NN\}
=\{T^k(x)\;:\;k\in\NN\}\;.
$$ 
The result then follows from the commutativity of the diagram
in Equation \eqref{commutTheta2}.

\medskip (2) {\bf Step 1 : } It is clear that the map
$(\Theta_2)_*:M_T\ra M_{\M}$ is order preserving and linear (for
positive scalars). Let us prove that it is a weak-star
homeomorphism. There is a direct proof of this, but we will use a
symbolic dynamics argument, as it is illuminating.

For every $x\in C$, let $[0;a_1,a_2,\dots]$ be the continued fraction
expansion of $x_+\in\M^\natural$ and let $[a_0;a_{-1},a_{-2},\dots]$
be the continued fraction expansion of $x_-\in({}^c\OOOO)^\natural$,
so that $a_i\in R-\FF_q$ for all $i\in\ZZ$. Let
$\sigma:(R-\FF_q)^\ZZ\ra (R-\FF_q)^\ZZ$ be the {\it (two-sided) shift}
defined by $\big(\sigma((x_i)_{i\in\ZZ})\big)_i=x_{i+1}$ for all
$i\in\ZZ$. Then by \cite[Theo.~3.7]{BroPau07JLMS} and Theorem
\ref{theo:commutfrmArtin} (2) and (5), the map $\Xi: C\ra
(R-\FF_q)^\ZZ$ defined by $x\mapsto (a_i)_{i\in\ZZ}$ is a
homeomorphism such that the following diagram commutes
$$
\begin{array}{ccc}
C &\stackrel{\Xi}{\longrightarrow}& (R-\FF_q)^\ZZ \medskip\\
T\downarrow\;\;\;\;& &\;\;\;\downarrow \sigma\medskip\\
C &\stackrel{\Xi}{\longrightarrow}& (R-\FF_q)^\ZZ\;.
\end{array}
$$ 
Let $\sigma_+:(R-\FF_q)^{\NN-\{0\}}\ra (R-\FF_q)^{\NN-\{0\}}$ be
the {\it (one-sided) shift} defined by the formula $\big(\sigma_+(
(x_i)_{i\in{\NN-\{0\}}})\big)_i=x_{i+1}$ for all $i\in\NN-\{0\}$. By
the properties of the continued fraction expansions, the map
$cfe:\M^\natural\ra (R-\FF_q)^{\NN-\{0\}}$ defined by $f\mapsto
(a_i(f))_{i\in\NN-\{0\}}$ is a homeomorphism such that the following
diagram commutes
$$
\begin{array}{ccc}
\M^\natural &\stackrel{cfe}{\longrightarrow}& (R-\FF_q)^{\NN-\{0\}} \medskip\\
\Psi\downarrow\;\;\;\;& &\;\;\;\downarrow \sigma_+\medskip\\
\M^\natural &\stackrel{cfe}{\longrightarrow}& (R-\FF_q)^{\NN-\{0\}}\;.
\end{array}
$$ 
Let $\pi_+:(R-\FF_q)^\ZZ\ra(R-\FF_q)^{\NN-\{0\}}$ be the map, called
the {\it natural extension}, forgetting the past, defined by
$\pi_+((x_i)_{i\in\ZZ})=(x_i)_{i\in{\NN-\{0\}}}$. The following
diagram commutes by construction
$$
\begin{array}{ccc}
C&\stackrel{\Xi}{\longrightarrow}& (R-\FF_q)^{\ZZ} \medskip\\
\Theta_2\downarrow\;\;\;\;& &\;\;\;\downarrow \pi_+\medskip\\
\M^\natural &\stackrel{cfe}{\longrightarrow}& (R-\FF_q)^{\NN-\{0\}}\;.
\end{array}
$$ 
Note that, as already seen, for every $x\in C$,
$$ 
t_1(x)=2\deg a_1(\Theta_2(x))\;,
$$ 
so that $\int_{\M} \deg a_1\;d(\Theta_2)_*\mu= \frac{1}{2} \int_{C}
t_1\,d\mu$ for every measure $\mu$ on $C$.

Step 1 now follows, by conjugation of dynamical systems, from the
well-known fact (see for instance \cite{Kitchens98}, or use the fact
that the symbolic Borel $\sigma$-algebras are generated by the
cylinders and Kolmogorov's extension theorem) that the pushforwards
map $(\pi_+)_*$ is a weak-star homeomorphism between the space of
$\sigma$-invariant measures on $(R-\FF_q)^{\ZZ}$ for which $\deg a_1$
is integrable and the space of $\sigma_+$-invariant measures on
$(R-\FF_q)^{\NN-\{0\}}$ for which $\deg a_1$ is integrable.

\medskip
\noindent 
{\bf Step 2 : } Let us prove that the map $F_2:M_T\ra
M_{T^\ZZ\uA(\OOOO)}$ is a weak-star homeomorphism.

Since $\int_{C\uA(\OOOO)} t_1\,d(\mu\otimes\mu_{\uA(\OOOO)})= \int_{C}
t_1\,d\mu$ for every measure $\mu$ on $C$, and since $T^\ZZ$ and
$\uA(\OOOO)$ commute, the measure $\mu\otimes\mu_{\uA(\OOOO)}$ is
indeed $T^\ZZ\times\uA(\OOOO)$-invariant if $\mu$ is $T$ invariant,
hence $F_2$ is well defined.  It is clear that $F_2$ is order
preserving and linear.

Every finite nonzero $T^\ZZ\times\uA(\OOOO)$-invariant measure $m$ on
$C\uA(\OOOO)$ disintegrates by the fibration $C\uA(\OOOO)\ra C$. The
conditional measures on the fibers, which are $\uA(\OOOO)$-invariant
and may be assumed to be probability measures, are hence equal to the
$\uA(\OOOO)$-invariant probability measures, by uniqueness. Since the
return time $t_1$ is constant on the fibers, the measure on $C$ over
which $m$ disintegrates hence belongs to $M_T$ if $m$ belongs to
$M_{T^\ZZ\uA(\OOOO)}$.  Step 2 follows.

\medskip
\noindent 
{\bf Step 3 : } Let us finally prove that the map $F_1:M_A\ra
M_{T^\ZZ\uA(\OOOO)}$ is a weak-star homeomorphism. Since it is clear
that $F_1$ is order preserving and linear, Steps 1 to 3 prove
Assertion (2) of Theorem \ref{theo:troispoint} by composition.

Since $C\uA(\OOOO)$ is open in $X^\natural$, and since any $A$-orbit
in $X^\natural$ meets $C$, for every finite nonzero $A$-invariant
measure $m$ on $X^\natural$, the measure $m_{\mid C\uA(\OOOO)}$ is
finite and nonzero. The open subset $C\uA(\OOOO)$ is a transversal for
the action of $\aaa_0$. If $t_0:X^\natural\ra\NN$ is the first passage
time in $C\uA(\OOOO)$ of the orbit of $\aaa_0$, that is,
$$
t_0:x\mapsto \min\{n\in\NN\;:\;\aaa_0^nx\in C\uA(\OOOO)\}\;,
\footnote{This map from $X^\natural$ to $\NN$ extends to $X^\natural$
  the map from $C\uA(\OOOO)$ to $\NN$ previously defined in Theorem
  \ref{theo:commutfrmArtin} (4) to be the zero map.}
$$ 
and if $T_0:X^\natural\ra C\uA(\OOOO)$ is the first passage map
$x\mapsto \aaa_0^{t_0(x)}x$ in $C\uA(\OOOO)$ of the orbit of $\aaa_0$,
then $\{t\in\ZZ\;:\; \aaa_0^tx\in C\uA(\OOOO)\}$ is exactly the set
$\{t_0(x)+t_k(T_0(x))\;:\; k\in\ZZ\}$. Since $C\uA(\OOOO)$ is open,
the disintegration of any $A$-invariant measure with respect to this
transversal is simply the restriction of $m$ to $C\uA(\OOOO)$, and
this restriction is $\uA(\OOOO)$-invariant. Therefore the discrete
time dynamical system $(X^\natural,\aaa_0)$ identifies with the
discrete time suspension of the discrete time dynamical system
$(C\uA(\OOOO),T)$ with roof function $t_1$. Furthermore, any measure
on $C\uA(\OOOO)$ which is invariant both under $T$ and $\uA(\OOOO)$
suspends to an $A$-invariant measure, whose total mass is the
integral of the roof function. Step 3 follows.

\medskip (3) The proof of the first claim in Assertion (3) of Theorem
\ref{theo:troispoint}, that is, the computation of 
$$
m_C = (F_2)^{-1}\circ F_1(m_X)\;,
$$
follows from the next proposition.

\bprop \label{prop:computmX}
The restriction $F_1(m_X)$ to $C\uA(\OOOO)$ of the $G$-invariant
probability measure $m_X$ on $X$ satisfies, for all $x\in C$ and $a\in
\uA(\OOOO)$, using the homeomorphism $xa\mapsto (x_-,x_+,a)$ from
$C\uA(\OOOO)$ to $({}^c\OOOO)^\natural\times \M^\natural\times
\uA(\OOOO)$,
$$
dm_X(xa)=\frac{q(q-1)^2}{2}\;\frac{d\mu_{\wh K}(x_-)}{|x_-|_\infty^4}\,
d\mu_{\wh K}(x_+)\,d\mu_{A}(a)\;.
$$
\eprop

\dem Let us denote by $dv$ the counting measure on the discrete set
$V\TT$. Recall that {\it Hopf's parametrisation} is the homeomorphism
from $\G\TT$ onto its image in $\partial_\infty
\TT\times\partial_\infty \TT\times V\TT$ defined by $\ell\mapsto
(\ell(-\infty), \ell(+\infty), \ell(0))$. Let $\wt m_{\rm BM}$ be the
(locally finite, positive, Borel) measure on $\G\TT$ defined, using
Hopf's parametrisation, on the full measure set of $\ell\in\G\TT$ such
that $\ell(\pm\infty)\neq\infty$, by
\begin{equation}\label{eq:bowenmargulis}
d\wt m_{\rm BM}(\ell)=\frac{d\mu_{\rm Hau}(\ell(-\infty))\;d\mu_{\rm
    Hau}(\ell(-\infty))
  \;dv(\ell(0))}{|\ell(-\infty)-\ell(+\infty)|^2_\infty}\;.
\end{equation}

By well-known arguments of for instance \cite{BurMoz96},
\cite{Roblin03}, \cite[\S 4]{BroPau07JLMS} and \cite[\S
  4.4]{BroParPau18}, since $\Ga$ is a lattice in $G$, the measure $\wt
m_{\rm BM}$ is invariant on the left by the action of $\Ga$ and on the
right by the geodesic flow, and defines a finite nonzero measure
$m_{\rm BM}$ on $\Ga\bs\G\TT$ invariant by the geodesic flow (called
the {\it Bowen-Margulis measure}). Furthermore, the $G$-invariant
probability measure $m_X$ on $X$ disintegrates with respect to the
fiber bundle $X=\Ga\bs G\ra \Ga\bs\G\TT=\Ga\bs G/\uA(\OOOO)$ defined
by $\Ga g\mapsto \Ga g\ell_*$, over the probability measure
$\frac{m_{\rm BM}}{\|m_{\rm BM}\|}$, with conditionnal measures the
$\uA(\OOOO)$-invariant probability measures on the fibers. 

Recall that for all $x\in C$, we have $x_-=g\cdot \infty=
g\ell_*(-\infty)$ and $x_+=g\cdot 0=g\ell_*(+\infty)$ where $g$ is the
unique lift of $x$ in $\wt C$; furthermore $g\ell_*(0)=g\,*=*$ is
constant. As seen during the proof of Theorem
\ref{theo:commutfrmArtin}, if $\ga\in\Ga$ fixes $g\ell_*$ for some
$g\in \wt C$, then it fixes $*$, $e_\infty$ and $e_0$, hence it
belongs to $\PGL_2(\OOOO)\cap\Ga=\PGL_2(\FF_q)$ and is diagonal. But a
nontrivial element of $\PGL_2(\FF_q)\cap A$ does not fix an
irrational point at infinity. Hence the stabilizer in $\Ga$ of an
element of $\wt C\,\ell_*$ is trivial. By Equation
\eqref{eq:bowenmargulis}, we thus have, for all $x\in C$ and $a\in
\uA(\OOOO)$,
$$
dm_{X}(xa) =\frac{1}{\|m_{\rm BM}\|}\;\frac{d\mu_{\rm Hau}(x_-)\;
d\mu_{\rm  Hau}(x_+)}{|x_--x_+|^2_\infty}\,d\mu_{A}(a)\;.
$$ 
For all $x\in C$, we have $x_-\in {}^c\OOOO$ and $x_+\in \M$, hence
$|x_+|_\infty <1< |x_-|_\infty$, so that $|x_--x_+|_\infty=
|x_-|_\infty$.  By \cite[Prop.~15.2 (1)]{BroParPau18} (where $\mu_{\rm
  Hau}$ was normalized in this reference to have total mass
$\frac{q+1}{q}$ instead of $1$, by [loc. cit., Prop.~15.2 (2)]), on
the set of $\xi\in\partial_\infty\TT=\wh K\cup\{\infty\}$ such that
$\xi\neq \infty$, we have
$$
d\mu_{\rm Hau}(\xi)=\frac{q}{q+1}\;
\frac{d\mu_{\wh K}(\xi)}{\max\{1,|\xi|_\infty^2\}}\;.
$$ 
Note that $\|m_{\rm BM}\|=\frac{2}{q(q-1)^2}\Big(\frac{q}{q+1}\Big)^2$,
by \cite[Prop.~15.3 (1)]{BroParPau18} (again using the above
difference of normalization).
Hence 
$$
dm_{X}(xa) =
\frac{q(q-1)^2}{2}\;\frac{d\mu_{\wh K}(x_-)}{|x_-|^4_\infty}\;
d\mu_{\wh K}(x_+)\,d\mu_{A}(a)\;.
$$ 
The result follows. \cqfd

\medskip By the first claim in Assertion (3) of Theorem
\ref{theo:troispoint}, we immediately have that there exists a
constant $c>0$ such that $(\Theta_2)_*m_C= c \;{\mu_{\wh K}}_{\;\mid
  \M^\natural} = c \;\mu_{\M}$. We have $\mu_{\wh K}(\M)=\frac{1}{q}$
by the normalization of $\mu_{\wh K}$ and the fact that
$\OOOO=\sqcup_{x\in\FF_q} (x+\M)$. This implies that $c=q\,\|m_C\|=
q\,m_X(C\uA(\OOOO))$.

Note that $\wh K-\OOOO=\sqcup_{u\in R-\FF_q}(u+\M)$,
and that there are $(q-1)q^n$ polynomials in $R$ with degree $n$ for
every $n\in\NN$.  Hence, by Proposition \ref{prop:computmX}, we have
\begin{align*}
m_X(C\uA(\OOOO))&=\frac{q(q-1)^2}{2}\;
\int_{{}^c\OOOO}\frac{d\mu_{\wh K}(x_-)}{|x_-|_\infty^4}\;
\mu_{\wh K}(\M)\;\mu_{A}(\uA(\OOOO))\\ &=\frac{(q-1)^2}{2}\;
\sum_{u\in R-\FF_q}\int_{y\in\M}\frac{d\mu_{\wh K}(y)}{|u+y|_\infty^4}
=\frac{(q-1)^2}{2}\;
\sum_{u\in R-\FF_q}\frac{\mu_{\wh K}(\M)}{|u|_\infty^4}\\ &=
\frac{(q-1)^2}{2q}\;
\sum_{n=1}^{+\infty}\frac{(q-1)q^n}{q^{4n}}=
\frac{q^3(q-1)^2}{2(q^2+q+1)}\;.
\end{align*}
This proves the second claim in Assertion (3) of Theorem
\ref{theo:troispoint}, and the last one follows by composition.

\medskip 
(4) Let us fix $x\in C$ which is $A$-periodic in $X^\natural$. As seen
in Assertion (1), $x$ is periodic under $T$. Up to replacing $x$ by
$T^kx$ for any $k$ big enough, since $T^kx\in xA$ so that $\mu_{T^kx}
=\mu_x$, we may assume by Equation \eqref{commutTheta2} that the
continued fraction expansion of $\Theta_2(x)$ is periodic. We then have
$$
C\uA(\OOOO)\cap xA=\{T^n(x)\uA(\OOOO)\;:\;0\leq n< n_x\}\;.
$$
Furthermore, for all $y\in C$ and $a\in\uA(\OOOO)$,
$$
d{\mu_{x}}_{\,\mid C\uA(\OOOO)}(ya) =\frac{\mu_x(C\uA(\OOOO))}{n_x}\;
\sum_{n=0}^{n_x-1}\delta_{T^n(x)}(y)\;d\mu_A(a)\;,
$$ 
where $\mu_x(C\uA(\OOOO))=\frac{n_x}{\sum_{n=0}^{n_x-1}
  t_1(T^n(x))}$ is the proportion of points of the $\aaa_0^\ZZ$-orbit
of $x$ that are in $C\uA(\OOOO)$.  Recalling that the first return
time satisfies $t_1(x)= 2\deg a_1(\Theta_2(x))$ and by the
commutativity of the diagram in Equation \eqref{commutTheta2}, this
proves the result.  \cqfd

\section{Behavior of continued fraction expansions along Hecke rays}
\label{sec:applyKPS}

In this section, we give more precisions on the correspondence between
$A$-orbits in the moduli space $X$ and $\Psi$-orbits in $\M$, for the
case of quadratic irrationals.  We study the behavior of the degrees
of the coefficients of the continued fraction expansions of the images
in $\M$ of $A$-periodic elements in $X$ varying along Hecke rays, and
we prove Theorem \ref{theo:escapingcfe} and Theorem
\ref{theo:intro:distribution}.

We start by recalling some definitions from for instance
\cite{KemPauSha17}. For every $x\in X$ and every prime polynomial
$P$ in $R$, the {\it $P$-Hecke tree} $T_P(x)$ with root $x$ is
the connected component of $x$ in the graph with vertex set $X$, with
an edge between the homothety classes of two $R$-lattices $\Lambda$
and $\Lambda'$ when $\Lambda'\subset \Lambda$ and $\Lambda/\Lambda'$
is isomorphic to $R/P R$ as an $R$-module. If $x=\Ga
g=g^{-1}[R\times R]$, we hence have
$$
T_P(x)=\{g^{-1}\ga\,[R\times P^n R]\;:\;\ga\in\Ga,\;n\in\NN\}\;.
$$ 
The boundary at infinity of $T_P(x)$ identifies with the projective
line $\PP_1(K_P)$ over the completion $K_P$ of $K$ for the $P$-adic
absolute value $|x|_P=q^{-v_P(x)}$, where $v_P(P^m\frac{A}{B})=m$ for
all $m\in\ZZ$ and $A,B\in R$ not divisible by $P$ (with $B\neq 0$). A
geodesic ray starting from $x$ (called a {\it Hecke ray}) in $T_P(x)$
is called {\it rational} if its point at infinity belongs to
$\PP_1(K)$, that is, if it is of the form $n\mapsto g^{-1}\ga[R\times
  P^n R]$ for some fixed $\ga\in\Ga$.
We refer for instance to \cite[\S 1]{KemPauSha17} for details and more
informations.

\medskip 
Now let us fix $f\in QI$.  We denote by $f^\sigma\in QI$ the Galois
conjugate of $f$ over $K$. Let $g_f= \begin{bmatrix} f^\sigma & f \\ 1
  & 1 \end{bmatrix} \in \PGL_2(\wh K)$. Let $t_f =
\frac{f^\sigma-1}{1-f}\in\wh K$ and $a_f=\alpha_0(t_f)
=\begin{bmatrix} 1 & 0 \\ 0 & t_f \end{bmatrix} \in A$.  Finally, let
$$
x_f= \Ga\, g_f\,a_f\in X\;.
$$ 
Note that for all $\begin{pmatrix} a & b \\ c & d \end{pmatrix}\in
\GL_2(K)$, we have
\begin{equation}\label{eq:calc2by2gf}
\begin{pmatrix} a & b \\ c & d\end{pmatrix}
\begin{pmatrix} f^\sigma & f \\ 1 & 1\end{pmatrix}
=\begin{pmatrix} \frac{af^\sigma+b}{cf^\sigma+d} &\frac{af+b}{cf+d}
 \\ 1 & 1\end{pmatrix} 
\begin{pmatrix} cf^\sigma+d & 0 \\ 0 & cf+d\end{pmatrix}\;.
\end{equation}
Hence for every $\ga\in \PGL_2(K)$, we have 
\begin{equation}\label{eq:calc2by2gfconseq}
x_{f} A =\Ga \ga\,g_{f}\,a_f A = \Ga g_{\ga\cdot f} A= x_{\ga\cdot f} A\;.
\end{equation}

\blemm \label{lem:xfAper} 
The points $\Ga\, g_f$ and $x_f$ of $X$ are $A$-periodic. If $f\in \M$
and $f^\sigma\in \,^c\OOOO$, then $a_f\in A_+$, $x_f\in C$ and
$\Theta(x_f)=(f^\sigma,f)$, so that $\Theta_2(x_f)=f$.  
\elemm

\dem Let $\ga \in\Ga-\{\id\}$ projectively fixing $f$ (hence
$f^\sigma$), with $v_\infty(\tr \wt\ga)\neq 0$ for some lift $\wt
\ga\in\GL_2(R)$ of $\ga$. The existence of $\ga$ is classical, see for
instance \cite{Serre83,KemPauSha17} and Proposition 17.2 with Equation
(15.6) in \cite{BroParPau18}.  By Equation \eqref{eq:calc2by2gf},
there exists $w\in R+fR$ such that $\ga \,g_f = g_f\begin{bmatrix}
w^\sigma & 0 \\ 0 & w \end{bmatrix}= g_f\,\alpha_0(t)$, with
$v_\infty(t) = v_\infty(\tr \wt\ga)\neq 0$. The fact that $\Ga\, g_f$
(and hence $x_f$, which is in the same $A$-orbit) is $A$-periodic then
follows for instance from \cite[Prop.~11]{KemPauSha17}.

Now assume that $f\in \M$ and $f^\sigma\in\,^c\OOOO$. Then
$v_\infty(t_f) = v_\infty(f^\sigma) <0$, hence $a_f\in A_+$.

Let $g'_f=g_f\,a_f$. An easy computation gives
$$
g'_f\cdot(\infty,0,1)=\begin{bmatrix} f^\sigma & f \\ 1 & 1 \end{bmatrix}
\begin{bmatrix} 1 & 0 \\ 0 & t_f \end{bmatrix}\cdot(\infty,0,1)=
(f^\sigma,f,1)\;.
$$ 
Hence $g'_f\in \wt C$, thus $x_f=\Ga g'_f \in C=\Ga \wt C$. We have
\begin{equation}\label{eq:surjfsigmaf}
\Theta(x_f)=\wt\Theta(g'_f)=
(g'_f\cdot\infty, \;g_f'\cdot 0)
=(g_f\cdot\infty,\;g_f\cdot 0)= (f^\sigma,f)\;,
\end{equation}
so that $\Theta_2(x_f)=\pr_2\circ\Theta(x_f)=f$.
\cqfd

\bigskip\noindent{\bf Proof of Theorem \ref{theo:escapingcfe}. } Let
$f\in QI$.  

Let us prove the first assertion of Theorem
\ref{theo:escapingcfe}.  We fix an irreducible polynomial $P\in R$.
For every $n\in\NN$, let $\ga_n\in\Ga$ be such that the quadratic
irrational
$$
f_n=\ga_n\cdot(P^n f)
$$ 
satisfies $f_n\in \M$ and ${f_n}^{\sigma}\in \,{}^c\OOOO$.  Since the
projective action of $\Ga=\PGL_2(R)$ on $QI$ does not change the
period of the eventually periodic continued fraction expansion of $f$,
up to cyclic permutation and multiplications by elements of
$\FF_q^\times$ (which do not change the degrees), there exist
$k_n,k'_n \in\NN-\{0\}$ such that the continued fraction expansions of
$f_n$ and $P^n f$ are indeed periodic after the times $k_n$ and $k'_n$
respectively, with the same period length $\ell_n\in\NN-\{0\}$, such
that for every $i\in\NN$,
$$
\deg a_{k_n+i}(P^nf)= \deg a_{k'_n+i}(f_n)\;.
$$

By Equation \eqref{eq:calc2by2gf}, for every $n\in\NN$, there exists
$b_n\in A$ such that $ g_{f_n} b_n= \ga_n \,g_{P^nf}$.  Hence
\begin{align*}
\Ga g_{f_n}b_n & =(g_{f_n}b_n)^{-1}[R\times R]=(\ga_n g_{P^nf})^{-1}[R\times R]=
g_{f}^{-1}\begin{bmatrix}
  P^n & 0 \\ 0 & 1 \end{bmatrix}^{-1}\ga_n^{-1}[R\times R]\\ & =
g_{f}^{-1}[R\times P^nR]\;.
\end{align*}
Therefore $(\Ga g_{f_n}b_n)_{n\in\NN}$ is the sequence of
vertices along a rational Hecke ray in the $P$-Hecke tree $T_P(\Ga
g_f)$. Note that $\Ga g_f$ is $A$-periodic by the first claim of Lemma
\ref{lem:xfAper}.

For every $n\in\NN$, the point $x_{f_n}=\Ga g_{f_n}a_{f_n}$, which
belongs to $C$ by Lemma \ref{lem:xfAper}, is in the same $A$-orbit as
$\Ga g_{f_n}b_n$. Let $\lambda_n\in\NN-\{0\}$ be the period of $\Ga
g_{f_n}\,a_{f_n}\,\ell_*$ under the geodesic flow on $\Ga\bs\G\TT$,
which is also the period of $\Ga g_{f_n}\,b_{n}\,\ell_*$, since $A$
preserves $\ell_*$.

Since $\ell_*(+\infty)=0$, the geodesic line $g_{f_n}\,a_{f_n}
\,\ell_*$, through $*$ at time $t=0$, has point at $+\infty$ equal,
using Lemma \ref{lem:xfAper}, to
$$
g_{f_n}\,a_{f_n}\cdot 0=\Theta_2(x_{f_n})=f_n.
$$ 
By Equation \eqref{eq:commuta0geodflow}, and by the geometric
interpretation of \cite[\S 6]{Paulin02}, we have
$$
\lambda_n=\sum_{i=k'_n}^{k'_n+\ell_n-1} \;2\;\deg \; a_i(f_n)=
\sum_{i=k_n}^{k_n+\ell_n-1} \;2\;\deg \; a_i(P^nf)\;.
$$ 
With $p_\infty: X \ra \Ga\bs V\TT$ the canonical projection given
by $\Ga g\mapsto \Ga g *$, we defined in \cite[Eq.~4]{KemPauSha17} the
{\it height} of $x\in X$ by
$$
\heit_\infty(x)= d_{\Ga\bs V\TT}\,(p_\infty(x),\Ga \,*)\;.
$$ 
Note that $\heit_\infty(x\,a)=\heit_\infty(x)$ for all $x\in X$ and
$a\in \uA(\OOOO)$. Hence using again
Equation \eqref{eq:commuta0geodflow} and the geometric interpretation
of \cite[\S 6]{Paulin02}, we have
\begin{align*}
\heit_\infty(\Ga g_{f_n}b_n)&\leq
\sup_{i\in\ZZ}\heit_\infty(\Ga g_{f_n}b_n\aaa_0^{i})
= \sup_{i\in\ZZ}\heit_\infty(\Ga g_{f_n}a_{f_n}\aaa_0^{i}) \\ & =
\max_{i=k'_n,\,\dots,\,k'_n+\ell_n-1} \;\deg \; a_i(f_n) =
\max_{i=k_n,\,\dots,\,k_n+\ell_n-1} \;\deg \; a_i(P^nf)\;.
\end{align*}

By \cite[Eq.~(14)]{KemPauSha17} and the four lines following it,
there exists $c>0$ such that for every $N\in\NN$,
$$
\liminf_{n\ra+\infty}
\frac{2\big(\max_{i=k_n,\,\dots,\,k_n+\ell_n-1} \;\deg \; a_i(P^nf)-N\big)}
{2\sum_{i=k_n}^{k_n+\ell_n-1} \;\deg \; a_i(P^nf)} \geq \liminf_{n\ra+\infty} 
\frac{2(\heit_\infty(\Ga g_{f_n}b_n)-N)}{\lambda_n}\geq c\;.
$$
This indeed proves the first assertion of Theorem 
\ref{theo:escapingcfe}.

\medskip
The second assertion of Theorem \ref{theo:escapingcfe} follows from
\cite[Lem.~13]{KemPauSha17}, whose hypothesis are exactly the ones of
this second assertion. The change between $\liminf$ in the definition
of $c$-degree-escaping continued fraction expansions and $\limsup$ in
the second assertion of Theorem \ref{theo:escapingcfe} comes from the
fact that we need in the final Remark of \cite[\S 4.1]{KemPauSha17} to
take a subsequence of the sequence of vertices along the rational
Hecke ray in order to obtain full escape of mass.

\medskip 
In order to prove the last assertion of Theorem \ref{theo:escapingcfe},
let again $P\in R$ be an irreducible polynomial. Since the $P$-Hecke
tree of $x=\Ga g\in X$ is 
$$
T_P(x)= \Big\{g^{-1}\ga^{-1}\alpha_0(P^n)\,[R\times R]\;:\;
\ga\in\Ga,\;n\in\NN\Big\}\;,
$$ 
since $\alpha_0(P^{-n})\ga g_f$ and $g_{P^{n}\ga\cdot f}$ have the
same $A$-orbit on the right by Equation \eqref{eq:calc2by2gfconseq},
Theorem \ref{theo:escapingcfe} (3) then follows from \cite[\S
  4.2]{KemPauSha17}.  \cqfd

\bigskip
\noindent{\bf Proof of Theorem \ref{theo:intro:distribution} in the
  introduction. }

We will actually prove the stronger statement that the uncountably
many sequences $(\ga''_n)_{n\in\NN}$ in $\Ga$ we are going to
construct satisfy the fact that for each one of them, there exists an
increasing sequence $(k_n)_{n\in\NN}$ in $\NN$ such that the images of
$\ga''_n$ and $\ga''_{n+1}$ in the group $\PGL_2(R/(P^{k_n}R))$, by
the reduction modulo $P^{k_n}$ of their coefficients, are equal.

Let $P\in R$ be an irreducible polynomial. Let $x\in X$ and $g\in G$
be such that $x=\Ga g$. We again refer to \cite[\S 2.5]{KemPauSha17}
for the basic facts we are using on the Hecke trees. For every point
at infinity $\xi$ of the Hecke tree $T_P(x)$, let
$(x_n^\xi)_{n\in\NN}$ be the sequence of vertices along the geodesic
ray in $T_P(x)$ from $x$ to $\xi$. In particular $x_n^0=g^{-1}[R\times
  P^nR]$. Recall that the action of $\Ga$ on the sphere of radius
$n\in\NN$ centered at $x$ in $T_P(x)$ is transitive. Hence there
exists a sequence $(\ga_n)_{n\in\NN}$ in $\Ga$ such that
$x_n^\xi=g^{-1}\ga_n^{-1} [R\times P^nR]$. The fact that $x_n^\xi$
belongs to the geodesic segment from $x$ to $x_{n+1}^\xi$ is
equivalent to the fact that $\ga_{n+1}\ga_n^{-1}$ belongs to the
mod$\; P^n$-upper triangular subgroup of $\Ga$, which is the subgroup
of elements of $\Ga$ whose $(2,1)$-coefficient vanishes modulo $P^n$.

By \cite[Theo.~4]{KemPauSha17}, there exists $c'>0$ and uncountably
many points at infinity $\xi$ of $T_P(x)$ such that for every
$A$-periodic $y\in X$, there exists an increasing sequence
$(n_k)_{k\in\NN}$ in $\NN$ such that the weak-star limit
$\theta=\lim_{k\ra\infty} \mu_{x_{n_k}}$ exists and satisfies
$\theta\geq c'\mu_y$.  Since $F:M_A\ra M_\M$ is weak-star continuous,
linear and order preserving by Theorem \ref{theo:troispoint} (2), we
have $\lim_{k\ra\infty} F(\mu_{x_{n_k}})\geq c'F(\mu_y)$. Hence since
the measures are finite and nonzero, with $c''=
\frac{c'\|F(\mu_y)\|}{\|\theta\|}>0$ (which depends on $\xi,y$ and the
subsequence), we have
\begin{equation}\label{eq:limmino}
\lim_{k\ra\infty} \frac{F(\mu_{x_{n_k}^\xi})}{\|F(\mu_{x_{n_k}^\xi}\|}
\geq c''\frac{F(\mu_y)}{\|F(\mu_y)\|}\;.
\end{equation}

Now let $f\in QI$.  We may assume that $f\in\M^\natural$ and
$f^\sigma\in ({}^c\OOOO)^\natural$, up to replacing $f$ by $\ga\cdot
f$ by some $\ga\in\Ga$, which does not change the conclusion of
Theorem \ref{theo:intro:distribution}. We now take $x=x_f\in C$ and
$g=g_f a_f$, so that $x_f=\Ga g_fa_f$ (see the lines above Lemma
\ref{lem:xfAper}). We consider the uncountably many sequences
$(\ga_n)_{n\in\NN}$ in $\Ga$ associated with the above uncountably
many points at infinity $\xi$ of $T_P(x_f)$.  Let $\beta_n\in\Ga$ be
such that 

$\bullet$~ $f_n=\beta_n\cdot(P^{n}\ga_n\cdot f)$ belongs to
$\M^\natural$ and ${f_n}^\sigma\in ({}^c\OOOO)^\natural$,

$\bullet$~ and besides, such that $\nu_{f_n}=\nu_{P^{n}\ga_n\cdot f}$
(in order to have this, we might need to multiply a $\beta_n$
satisfying the first point by an element of $\alpha_0
(\FF_q^\times)\subset \Ga$, see \cite[Theo.~1]{BerNak00} for an
explanation).

Since $\beta_n\alpha_0(P^{-n})\ga_n\in\PGL_2(K)$ and by Equation
\eqref{eq:calc2by2gf}, there exists $b_n\in A$ such that
$$
\beta_n\alpha_0(P^{-n})\ga_n g_f=g_{\beta_n\alpha_0(P^{-n})\ga_n\cdot f} b_n 
=g_{f_n}b_n\;.
$$
Hence
$$
x_n^\xi=(g_f a_f)^{-1}\ga_n^{-1}[R\times P^nR]=
\Ga \alpha_0(P^{-n})\ga_ng_f a_f=\Ga g_{f_n}b_na_f\;,
$$
which is in the same $A$-orbit as $\Ga g_{f_n}a_{f_n}$. By Lemma
\ref{lem:xfAper}, we have $\Ga g_{f_n} a_{f_n}\in C$ and $\Theta_2(\Ga
g_{f_n} a_{f_n})=f_n$.  Hence by Theorem \ref{theo:troispoint} (4), we
have
$$
\frac{F(\mu_{x_n^\xi})}{\|F(\mu_{x_n^\xi}\|}=
\frac{F(\mu_{\Ga g_{f_n} a_{f_n}})}{\|F(\mu_{\Ga g_{f_n} a_{f_n}}\|}=
\nu_{f_n}=\nu_{P^{n}\ga_n\cdot f}\;.
$$

Let $f'\in QI$ and $y=x_{f'}$, which is $A$-periodic. Up to replacing
$f'$ by $\ga\cdot f'$ by some $\ga\in\Ga$, we assume that $f'\in\M$
and $(f')^\sigma\in ({}^c\OOOO)^\natural$, so that $\Theta_2(y)=f'$
again by Lemma \ref{lem:xfAper}.  Similarly, by Theorem
\ref{theo:troispoint} (4), we have
$$
\frac{F(\mu_{y})}{\|F(\mu_{y}\|}=\nu_{f'}\;.
$$

The result then follows from Equation \eqref{eq:limmino}.
\cqfd

%
%
%

\cqfd {\small \bibliography{../biblio} }

\bigskip
{\small
\noindent \begin{tabular}{l}
Laboratoire de math\'ematique d'Orsay,
\\ UMR 8628 Univ. Paris-Sud, CNRS\\
Universit\'e Paris-Saclay,
91405 ORSAY Cedex, FRANCE\\
{\it e-mail: frederic.paulin@math.u-psud.fr}
\end{tabular}
\medskip

\noindent \begin{tabular}{l} 
Mathematics Department, Technion \\
Israel Institute of Technology, Haifa, 32000 ISRAEL.\\
{\it e-mail: ushapira@tx.technion.ac.il}
\end{tabular}
}

\end{document}